\documentclass[10pt]{amsart}
\usepackage{amsfonts}
\usepackage{ifthen}
\usepackage{amsthm}
\usepackage{amsmath}
\usepackage{graphicx}
\usepackage{amscd,amssymb,amsthm}

\newcounter{minutes}
\divide\time by 60
\newcounter{hours}
\multiply\time by 60 \addtocounter{minutes}{-\time}

\newtheorem{theorem}{Theorem}
\newtheorem{lemma}{Lemma}

\keywords{Bessel functions of the first kind; modified Bessel functions of the first kind; close-to-convex functions; starlike functions; transcendental entire functions; zeros of cross product of Bessel functions; infinite product.} \subjclass[2010]{33C10, 30C45.}

\title[]{Starlikeness of a cross-product of Bessel functions}

\author[H.A. Al-kharsani]{Huda A. Al-Kharsani}
\address{Department of Mathematics, Girls College, University of Dammam, 838 Dammam, Saudi Arabia}
\email{halkharsani@ud.edu.sa}

\author[\'A. Baricz]{\'Arp\'ad Baricz}
\address{Department of Economics, Babe\c{s}-Bolyai University, Cluj-Napoca 400591, Romania}
\address{Institute of Applied Mathematics, \'Obuda University, 1034 Budapest, Hungary}
\email{bariczocsi@yahoo.com}

\author[T.K. Pog\'any]{Tibor K. Pog\'any}
\address{Faculty of Maritime Studies, University of Rijeka, 51000 Rijeka, Croatia}
\address{Institute of Applied Mathematics, \'Obuda University, 1034 Budapest, Hungary}
\email{poganj@pfri.hr}

\thanks{$^{\bigstar}$The research of H.A. Al-kharsani was supported by Deanship of Scientific Research in University of Dammam, project ID$\#$2014219. The work of \'A. Baricz was supported by a research grant of the Romanian National Authority for Scientific Research, CNCS-UEFISCDI, project number PN-II-RU-TE-2012-3-0190. The research of T.K. Pog\'any was covered in part by the Croatian
Science Foundation under the project No. 5435.}

\begin{document}

\def\thefootnote{}
\footnotetext{ \texttt{File:~\jobname .tex,
          printed: \number\year-\number\month-\number\day,
          \thehours.\ifnum\theminutes<10{0}\fi\theminutes}
} \makeatletter\def\thefootnote{\@arabic\c@footnote}\makeatother

\maketitle

\begin{center}
Dedicated to the memory of Professor Lee Lorch
\end{center}

\begin{abstract}
In this paper a necessary and sufficient condition is deduced for the close-to-convexity of a cross product of Bessel and modified Bessel functions of the first kind and their derivatives by using a result of Shah and Trimble about transcendental entire functions with univalent derivatives, the newly discovered power series and infinite product representation of this cross-product, as well as a slightly modified version of a result of Lorch on the monotonicity of the zeros of the cross product with respect to the order.
\end{abstract}

\section{\bf Introduction and the Main Results}

Let $J_{\nu}$ and $I_{\nu}$ denote the Bessel and modified Bessel functions of the first kind. Motivated by their
appearance as eigenvalues in the clamped plate problem for the
ball, Ashbaugh and Benguria have conjectured
that the positive zeros of the function $z\longmapsto J_{\nu}(z)I_{\nu}'(z)-J_{\nu}'(z)I_{\nu}(z)$
increase with $\nu$ on $\left[-\frac{1}{2},\infty\right).$ Lorch \cite{lorch} verified this conjecture and presented some other properties of
the zeros of the above cross product of Bessel and modified Bessel functions. In this paper we point out that actually the above monotonicity property is valid on $(-1,\infty).$ We are also interested on an application of the above monotonicity property of the zeros of the cross product of Bessel and modified Bessel functions. Namely, our aim is to find a necessary and sufficient condition on the parameter $\nu$ such that the normalized form of the above cross product maps the open unit disk into a starlike domain and all of its derivatives are close-to-convex, and hence univalent. It is worth to mention that geometric properties, like univalence, starlikeness, spirallikeness and convexity of Bessel functions were studied in the sixties by Brown \cite{brown}, and also by Kreyszig and Todd \cite{todd}. Note that some other geometric properties of Bessel functions of the first kind were studied later by others, see for example the papers \cite{publ,lecture,bsk,samy,barsza,szasz,szasz2} and the references therein. In order to prove our main results we use a result of Shah and Trimble \cite{st} about transcendental entire functions with univalent derivatives and power series and infinite product representations for the above cross product of Bessel and modified Bessel functions of the first kind. Our first main result is the following theorem.

\begin{theorem}\label{th1}
The function
$$z\longmapsto 2^{2\nu}z^{-\frac{\nu}{2}+\frac{3}{4}}\Gamma(\nu+1)\Gamma(\nu+2)\left(J_{\nu+1}(\sqrt[4]{z}) I_{\nu}(\sqrt[4]{z}) + J_\nu(\sqrt[4]{z})I_{\nu+1}(\sqrt[4]{z})\right)$$
is starlike in $\mathbb{D}$ and all of its derivatives are
close-to-convex (and hence univalent) there if and only if $\nu\geq \nu^{\ast},$ where $\nu^{\ast}\simeq-0.9427\dots$ is the unique root of the next equation on $(-1,\infty)$
$$(\nu-1)J_{\nu}(1)I_{\nu+1}(1)+(\nu-1)J_{\nu+1}(1)I_{\nu}(1)=J_{\nu}(1)I_{\nu}(1).$$
\end{theorem}

Motivated by the above result we also consider simply the product $J_{\nu}(z)I_{\nu}(z)$ in the next theorem.

\begin{theorem}\label{th2}
The function
$$z\longmapsto 2^{2\nu}z^{-\frac{\nu}{2}+1}\Gamma^2(\nu+1)J_{\nu}(\sqrt[4]{z})I_{\nu}(\sqrt[4]{z})$$
is starlike in $\mathbb{D}$ and all of its derivatives are
close-to-convex (and hence univalent) there if and only if $\nu\geq \nu^{\star},$ where $\nu^{\star}\simeq-0.4336\dots$
is the unique root of the next equation on $(-1,\infty)$
$$J_{\nu+1}(1)I_{\nu}(1)-J_{\nu}(1)I_{\nu+1}(1)=(\nu+1)J_{\nu}(1)I_{\nu}(1).$$
\end{theorem}

We note that if consider the particular cases
$$J_{-\frac{1}{2}}(z)=\sqrt{\frac{2}{\pi}}\frac{1}{\sqrt{z}}\cos z, \ J_{\frac{1}{2}}(z)=\sqrt{\frac{2}{\pi}}\frac{1}{\sqrt{z}}\sin z$$
and $$I_{-\frac{1}{2}}(z)=\sqrt{\frac{2}{\pi}}\frac{1}{\sqrt{z}}\cosh z,\ I_{\frac{1}{2}}(z)=\sqrt{\frac{2}{\pi}}\frac{1}{\sqrt{z}}\sinh z,$$
then by using the above theorems in particular for $\nu=-\frac{1}{2}$ we obtain that the functions
$$z\longmapsto \frac{1}{2}z^{\frac{3}{4}}\left(\sin\sqrt[4]{z}\cosh\sqrt[4]{z}+\cos\sqrt[4]{z}\sinh\sqrt[4]{z}\right)\ \ \mbox{and}\ \
z\longmapsto z\cos\sqrt[4]{z}\cosh\sqrt[4]{z}$$
are starlike in $\mathbb{D}$ and all of its derivatives are
close-to-convex (and hence univalent) there.

In order to prove our main results we will need some preliminary results. We note that although Lemma \ref{lem2}, \ref{lem3} and \ref{lem4} were deduced to prove Theorem \ref{th1}, they are also of independent interest and can be applied to solve other problems related to Bessel and modified Bessel functions of the first kind.

The next result of Shah and Trimble \cite[Theorem 2]{st} is one of the key tools in the proof of the main results.

\begin{lemma}\label{lem1}
Let $\mathbb{D}=\{z\in\mathbb{C}:|z|<1\}$ be the unit disk and $f:\mathbb{D}\to\mathbb{C}$ be an entire function of the form
$$f(z)=z\prod_{n\geq 1}\left(1-\frac{z}{z_n}\right),$$
where all $z_n$ have the same argument and satisfy $|z_n|>1.$ Then $f$ is starlike in $\mathbb{D}$ and all of its derivatives are
close-to-convex there if and only if the following inequality is valid
$$\sum_{n\geq1}\frac{1}{|z_n|-1}\leq 1.$$
\end{lemma}

It is worth to mention that by using the known recurrence relations $zJ_{\nu}'(z)-\nu J_{\nu}(z)=-zJ_{\nu+1}(z)$ and $zI_{\nu}'(z)-\nu I_{\nu}(z)=zI_{\nu+1}(z)$ the cross product $J_{\nu}(z)I_{\nu}'(z)-J_{\nu}'(z)I_{\nu}(z)$ actually can be rewritten as $J_{\nu+1}(z) I_\nu(z) + J_\nu(z)I_{\nu+1}(z).$ In the sequel we will also use the following power series and infinite product representations of the above cross product of Bessel and modified Bessel functions. These results complement the well-known results on Bessel and modified Bessel functions of the first kind and may be of independent interest.

\begin{lemma}\label{lem2}
If $\nu>-1$ and $z\in\mathbb{C},$ then we have the next power series representation
$$J_{\nu+1}(z) I_{\nu}(z) + J_{\nu}(z)I_{\nu+1}(z)=2\sum_{n \geq 0}\frac{(-1)^n\left( \frac{z}{2}\right)^{2\nu+4n+1}}{n!\Gamma(\nu+n+1)\Gamma(\nu+2n+2)}.$$
\end{lemma}

\begin{lemma}\label{lem3}
If $\nu>-1$ and $z\in\mathbb{C},$ then we have the next Hadamard factorization
$$2^{2\nu}z^{-2\nu-1}\Gamma(\nu+1)\Gamma(\nu+2)\left(J_{\nu+1}(z) I_{\nu}(z) + J_{\nu}(z)I_{\nu+1}(z)\right)=\prod_{n\geq 1}\left(1-\frac{z^4}{\gamma_{\nu,n}^4}\right),$$
where $\gamma_{\nu,n}$ is the $n$th positive zero of the function $z\longmapsto J_{\nu+1}(z) I_{\nu}(z) + J_\nu(z)I_{\nu+1}(z).$ Moreover, the zeros $\gamma_{\nu,n}$ satisfy the interlacing inequalities $j_{\nu,n}<\gamma_{\nu,n}<j_{\nu,n+1}$ and $j_{\nu,n}<\gamma_{\nu,n}<j_{\nu+1,n}$ for $n\in\mathbb{N}$ and $\nu>-1,$
where $j_{\nu,n}$ stands for the $n$th positive zero of the Bessel function $J_{\nu}.$
\end{lemma}

It is worth to mention that the inequalities $j_{\nu,n}<\gamma_{\nu,n}<j_{\nu+1,n}$ were proved for $\nu\geq -\frac{1}{2}$ by Lorch \cite{lorch}. To prove our main result in Theorem \ref{th1} we will need also the following result. Note that this result was proved also earlier by Lorch \cite{lorch} for the case $\nu\geq -\frac{1}{2}.$ Our proof is just a slight modification of the proof made by Lorch \cite{lorch}.

\begin{lemma}\label{lem4}
The positive zeros of $z\longmapsto J_{\nu}(z)I_{\nu}'(z)-J_{\nu}'(z)I_{\nu}(z)$ increase with $\nu$ on $(-1,\infty).$
\end{lemma}

\section{\bf Proofs of the Preliminary and Main Results}
\setcounter{equation}{0}

In this section our aim is to present the proof of the preliminary and main results.

\begin{proof}[\bf Proof of Lemma \ref{lem2}]
Let us consider the product of Bessel and modified Bessel functions of the first kind
   \[ \Pi_{\mu, \nu}(z) = \left( \frac2z\right)^{\mu+\nu}\, J_\mu(z) I_\nu(z),\]
where $\mu,\nu>-1$ and $z\in\mathbb{C}.$ We start with the series representations for both Bessel functions in the manner of Watson \cite[p. 147]{watson}. Thus
   \begin{equation} \label{A-1}
	    \Pi_{\mu, \nu}(z)  = \frac1{\Gamma(\nu+1)} \sum_{n \geq 0}
	                        \frac{(-1)^n\left( \frac z2\right)^{2n}}{n!\Gamma(\mu+n+1)}{}_2F_1(-n, -\mu -n; \nu+1; -1),
	 \end{equation}
where ${}_2F_1(a,b;c;\cdot)$ stands for the Gaussian hypergeometric function. Indeed, we have
   \begin{align*}
	    J_\mu(x) I_\nu(z) &= \left( \frac z2\right)^{\mu+\nu} \sum_{m \geq 0} \sum_{k \geq 0}\frac{(-1)^m\left( \tfrac z2\right)^{2(m+k)}}
			                     {m!k!\Gamma(\mu+m+1)\Gamma(\nu+k+1)} \\
												&= \left( \frac z2\right)^{\mu+\nu} \sum_{n \geq 0} \left\{ \sum_{m=0}^n  \frac{(-1)^m}
												   {(m-n)!m!\Gamma(\mu+m-n+1)\Gamma(\nu+m+1)}\right\}\left( \frac z2\right)^{2n} \\
												&= \left( \frac z2\right)^{\mu+\nu} \sum_{n \geq 0} \frac1{n!} \left\{ \sum_{m=0}^n \binom{n}{m} \frac{(-1)^m}
												   {\Gamma(\mu+n-m+1)\Gamma(\nu+m+1)}\right\}\left( \frac z2\right)^{2n}\, .
	 \end{align*}
As to the proof of \eqref{A-1}, it follows immediately from
   \[ \sum_{m=0}^n \binom{n}{m} \frac{(-1)^m}{\Gamma(\mu+n-m+1)\,\Gamma(\nu+m+1)}
	           = \frac{{}_2F_1(-n, -\mu -n; \nu+1; -1)}{\Gamma(\nu+1)\Gamma(\mu+n+1)}\,. \]				
Now, consider the Jacobi polynomial (or hypergeometric polynomial) \cite[p. 99]{askey}
   \[ P_n^{(\alpha, \beta)}(z) = \begin{cases}
	          \dfrac{(1+\alpha)_n}{n!}{}_2F_1\left(-n, n+\alpha+\beta+1; \alpha+1; \dfrac{1-z}2 \right) \\
						\dfrac{(1+\alpha)_n}{n!}\left( \dfrac{1+z}2\right)^n {}_2F_1\left(-n, -n-\beta; \alpha+1; \dfrac{z-1}{z+1} \right)
				\end{cases}\, ,	\]
where the latter expression \cite[p. 779]{AbramSte} is obtained from the first definition via the Pfaff transform of the hypergeometric term. This in turn implies that
$$\Pi_{\mu, \nu}(z) = \sum_{n \geq 0} \frac{(-1)^n P_n^{(\nu, \mu)}(0)}{2^n\Gamma(\mu+n+1) \Gamma(\nu+n+1)}z^{2n}.$$
Now, we are looking for a closed form in the case of the symmetric sum
   \begin{equation} \label{A1}
      \Pi_{\nu+1, \nu}(z) + \Pi_{\nu, \nu+1}(z) = \sum_{n \geq 0}
	                               \frac{(-1)^n\left[ P_n^{(\nu, \nu+1)}(0)+ P_n^{(\nu+1, \nu)}(0)\right]}
																 {2^n\Gamma(\nu+n+1) \Gamma(\nu+n+2)}z^{2n} .
	 \end{equation}
By using the recurrence relation \cite[p. 782]{AbramSte}
   \[ (1-z)P_n^{(\alpha+1, \beta)}(z) + (1+z)P_n^{(\alpha, \beta+1)}(z) = 2P_n^{(\alpha, \beta)}(z)\]
and the initial value
   \[ P_n^{(\nu, \nu)}(0) = \frac{\sqrt{\pi}\, \Gamma(\nu+n+1)}{\Gamma\left( \tfrac12 - \tfrac n2\right)\, n!\,
	                          \Gamma(\nu + \tfrac n2+1)}\, ,\]
together with
   \[ \frac1{\Gamma\left( \tfrac12 - \tfrac n2\right)} = \begin{cases}
	                 \dfrac{(-1)^k\, \prod_{j=1}^k(2j-1)}{\sqrt{\pi}\, 2^k}, & \mbox{if}\ \ n=2k,\, k\in \mathbb N \\
									 0,  & \mbox{if}\ \ n = 2k-1,\, k\in \mathbb N
							\end{cases},\]
for all $k\in \mathbb N_0$ we conclude that
   \[ P_n^{(\nu+1, \nu)}(0) +  P_n^{(\nu, \nu+1)}(0) = 2\, P_n^{(\nu, \nu)}(0) = \begin{cases}
	                 \dfrac{2\,\sqrt{\pi}\,(-1)^k\,(2k-1)!!\, \Gamma(\nu+2k+1)}{(2k)!\, 2^k\,\Gamma\big(\nu+ k+1\big)}, & \mbox{if}\ \ n = 2k\\
									 0, & \mbox{if}\ \ n= 2k-1
							\end{cases}\,.\]
This in conjunction with \eqref{A1} gives
   \begin{align*}
	    \Pi_{\nu+1, \nu}(z) + \Pi_{\nu, \nu+1}(z) &=  2 \sum_{n \geq 0}
	                                 \frac{(-1)^n \,(2n-1)!!\Gamma(\nu+2n+1)z^{4n}}{2^{3n}\, \Gamma(\nu+2n+1) \Gamma(\nu+2n+2)
													         \, (2n)!\, \Gamma\left(\nu+ n+1\right)} \nonumber \\
																&= 2 \sum_{n \geq 0} \frac{(-1)^n }{ \Gamma(\nu+n+1) \Gamma(\nu+2n+2) \, n!}
																   \, \left( \frac z2\right)^{4n}\, .
	 \end{align*}
Thus we have
   \[ J_{\nu+1}(z) I_\nu(z) + J_\nu(z)I_{\nu+1}(z) = 2\left( \frac{z}{2}\right)^{2\nu+1} \sum_{n \geq 0}
	                        \frac{(-1)^n\left(\frac{z}{2}\right)^{4n} }{n!\Gamma(\nu+n+1) \Gamma(\nu+2n+2)}. \]\end{proof}

\begin{proof}[\bf Proof of Lemma \ref{lem3}]
By following Lorch's approach \cite{lorch} we first show that for $\nu>-1$ the zeros $\gamma_{\nu,n}$ exist. By using the discussion in the introduction after Lemma \ref{lem1} it is clear that $\gamma_{\nu,n}$ is also the $n$th positive zero of the function
$$\varphi_{\nu}(z)=\frac{J_{\nu+1}(z)}{J_{\nu}(z)}+\frac{I_{\nu+1}(z)}{I_{\nu}(z)}.$$
In view of the Mittag-Leffler expansions
$$\frac{J_{\nu+1}(z)}{J_{\nu}(z)}=\sum_{n\geq 1}\frac{2z}{j_{\nu,n}^2-z^2},\ \ \frac{I_{\nu+1}(z)}{I_{\nu}(z)}=\sum_{n\geq 1}\frac{2z}{j_{\nu,n}^2+z^2}$$
it follows that
$$\varphi_{\nu}'(z)=\sum_{n\geq 1}\frac{4j_{\nu,n}^2(3j_{\nu,n}^4+z^4)}{(j_{\nu,n}^4-z^4)^2}>0$$
for $z\in\Delta=(0,j_{\nu,1})\cup (j_{\nu,1},j_{\nu,2})\cup{\dots}\cup(j_{\nu,n},j_{\nu,n+1})\cup{\dots}$ and $\nu>-1.$ On the other hand we have the following limits $\lim_{z\searrow0}\varphi_{\nu}(z)=0,$ $\lim_{z\nearrow j_{\nu,1}}\varphi_{\nu}(z)=\infty,$ $\lim_{z\searrow j_{\nu,1}}\varphi_{\nu}(z)=-\infty,$ $\lim_{z\nearrow j_{\nu,2}}\varphi_{\nu}(z)=\infty,\ {\dots},$ which together with the above monotonicity property implies that $j_{\nu,n}<\gamma_{\nu,n}<j_{\nu,n+1}$ for $n\in\mathbb{N}$ and $\nu>-1,$ that is, the zeros $\gamma_{\nu,n}$ and $j_{\nu,n}$ interlace. Moreover, since $\varphi_{\nu}(j_{\nu+1,n})=I_{\nu+1}(j_{\nu+1,n})/I_{\nu}(j_{\nu+1,n})>0$ for $n\in\mathbb{N}$ and $\nu>-1,$ it follows that $\gamma_{\nu,n}<j_{\nu+1,n}$ for $n\in\mathbb{N}$ and $\nu>-1.$ With these we proved the existence of the zeros and also their bounds.

Now, let us focus on the infinite product. We will show that for $\nu>-1$ and $z\in\mathbb{C}$ we have
\begin{equation}\label{prod}2^{2\nu}z^{-\nu-\frac{1}{2}}\Gamma(\nu+1)\Gamma(\nu+2)\left(J_{\nu+1}(\sqrt{z}) I_\nu(\sqrt{z}) + J_\nu(\sqrt{z})I_{\nu+1}(\sqrt{z})\right)=\prod_{n\geq 1}\left(1-\frac{z^2}{\gamma_{\nu,n}^4}\right).\end{equation}
This is equivalent to the first statement of Lemma \ref{lem3}. In view of Lemma \ref{lem2} we have
$$\frac{2^{2\nu}}{z^{\nu+\frac{1}{2}}}\Gamma(\nu+1)\Gamma(\nu+2)\left(J_{\nu+1}(\sqrt{z}) I_\nu(\sqrt{z}) + J_\nu(\sqrt{z})I_{\nu+1}(\sqrt{z})\right)=
1+\sum_{n \geq 1}\frac{(-1)^n\Gamma(\nu+1)\Gamma(\nu+2){z}^{2n}}{n!2^{4n}\Gamma(\nu+n+1)\Gamma(\nu+2n+2)}.$$
This is an entire function of growth order $\frac{1}{4}$ since
$$\lim_{n\to\infty}\frac{n\log n}{\log\Gamma(n+1)+\log\Gamma(\nu+n+1)+\log\Gamma(\nu+2n+2)+\log\frac{2^{4n}}{\Gamma(\nu+1)\Gamma(\nu+2)}}=\frac{1}{4}.$$
Here we used the limit $\frac{\log\Gamma(an+b)}{n\log n}\to a,$ as $n\to\infty,$ where $a,b>0.$ To see this just observe that
$$\lim_{x\to\infty}\frac{\log\Gamma(ax+b)}{x\log x}=a\lim_{x\to\infty}\frac{\psi(ax+b)}{1+\log x}=a\lim_{x\to\infty}\frac{\log(ax+b)-\frac{1}{2(ax+b)}+ \mathcal{O}(x^{-2})}{\log x}=a,$$
where $\psi(x)=\Gamma'(x)/\Gamma(x)$ is the logarithmic derivative of the Euler gamma function. Now, applying the Hadamard theorem \cite[p. 26]{lev} it follows that \eqref{prod} is indeed valid.
\end{proof}

\begin{proof}[\bf Proof of Lemma \ref{lem4}]
We know that the function $\nu\mapsto \varphi_{\nu}(z)$ is decreasing on $(-1,\infty)$ for $z>0$ fixed, see \cite{lorch}. Note that this can be verified also by looking at the Mittag-Leffler expansions in the above proof and by using the fact that the zeros $j_{\nu,n}$ are increasing on $(-1,\infty)$ as $\nu$ increases for each $n\in\mathbb{N}$ fixed. Hence, for $\varepsilon>0$ and $\nu>-1$ we have that $\varphi_{\nu+\varepsilon}(\gamma_{\nu,n})<\varphi_{\nu}(\gamma_{\nu,n})=0$ for some fixed $n\in\mathbb{N}.$ On the other hand, from the proof of Lemma \ref{lem3} we know that $\varphi_{\nu}$ is increasing on $\Delta$ for $\nu>-1$ and consequently we have that the equation $\varphi_{\nu+\varepsilon}(\gamma_{\nu+\varepsilon,n})=0$ together with above inequality yield $\gamma_{\nu+\varepsilon,n}>\gamma_{\nu,n},$ which completes the proof.
\end{proof}

\begin{proof}[\bf Proof of Theorem \ref{th1}]
By using Lemma \ref{lem3} we get that
$$2^{2\nu}z^{-\frac{\nu}{2}+\frac{3}{4}}\Gamma(\nu+1)\Gamma(\nu+2)\left(J_{\nu+1}(\sqrt[4]{z}) I_{\nu}(\sqrt[4]{z}) + J_\nu(\sqrt[4]{z})I_{\nu+1}(\sqrt[4]{z})\right)=z\prod_{n\geq 1}\left(1-\frac{z}{\gamma_{\nu,n}^4}\right).$$
On the other hand, by taking the logarithmic derivative of both sides of the relation
$$2^{2\nu}z^{-2\nu-1}\Gamma(\nu+1)\Gamma(\nu+2)\Phi_{\nu}(z)=\prod_{n\geq 1}\left(1-\frac{z^4}{\gamma_{\nu,n}^4}\right),$$
where $$\Phi_{\nu}(z)=J_{\nu+1}(z) I_\nu(z) + J_\nu(z)I_{\nu+1}(z),$$ we obtain that
$$\frac{1}{4}\left(2\nu+1-\frac{z\Phi_{\nu}'(z)}{\Phi_{\nu}(z)}\right)=\sum_{n\geq 1}\frac{z^4}{\gamma_{\nu,n}^4-z^4}.$$
Now, by using the well-known recurrence relations $zJ_{\nu+1}'(z)=zJ_{\nu}(z)-(\nu+1)J_{\nu+1}(z),$ $zI_{\nu+1}'(z)=zI_{\nu}(z)-(\nu+1)I_{\nu+1}(z),$ $zJ_{\nu}'(z)=\nu J_{\nu}(z)-zJ_{\nu+1}(z)$ and $zI_{\nu}'(z)=\nu I_{\nu}(z)+zI_{\nu+1}(z),$ it follows that
$$\Phi_{\nu}'(z)=2J_{\nu}(z)I_{\nu}(z)-\frac{1}{z}J_{\nu}(z)I_{\nu+1}(z)-\frac{1}{z}J_{\nu+1}(z)I_{\nu}(z).$$
Since in view of Lemma \ref{lem4} the function $\nu\mapsto \gamma_{\nu,n}$ is increasing on $(-1,\infty)$ for each $n\in\mathbb{N},$ it follows that for $n\in\{2,3,\dots\}$ and $\nu\geq\nu^{\ast}$ we have $\gamma_{\nu,n}>{\dots}>\gamma_{\nu,1}\geq \gamma_{\nu^{\ast},1}\simeq 1.1639{\dots}>1.$ Moreover, the above monotonicity property of the zeros $\gamma_{\nu,n}$ implies that
$$\nu\longmapsto\sum_{n\geq 1}\frac{1}{\gamma_{\nu,n}^4-1}$$
is decreasing on $(-1,\infty)$ and consequently
$$\sum_{n\geq 1}\frac{1}{\gamma_{\nu,n}^4-1}\leq1$$
if and only if $\nu\geq \nu^{\ast},$ where $\nu^{\ast}$ is the unique root of the equation
$$\sum_{n\geq 1}\frac{1}{\gamma_{\nu,n}^4-1}\leq1\ \ \Longleftrightarrow \ \ (\nu-1)J_{\nu}(1)I_{\nu+1}(1)+(\nu-1)J_{\nu+1}(1)I_{\nu}(1)=J_{\nu}(1)I_{\nu}(1).$$
Thus, applying Lemma \ref{lem1} the proof of this theorem is complete.
\end{proof}

\begin{proof}[\bf Proof of Theorem \ref{th2}]
In view of the well-known infinite products
$$2^{\nu}\Gamma(\nu+1)z^{-\nu}J_{\nu}(z)=\prod_{n\geq 1}\left(1-\frac{z^2}{j_{\nu,n}^2}\right), \ \
2^{\nu}\Gamma(\nu+1)z^{-\nu}I_{\nu}(z)=\prod_{n\geq 1}\left(1+\frac{z^2}{j_{\nu,n}^2}\right)$$
it follows that
$$2^{2\nu}z^{-\frac{\nu}{2}+1}\Gamma^2(\nu+1)J_{\nu}(\sqrt[4]{z})I_{\nu}(\sqrt[4]{z})=z\prod_{n\geq1}\left(1-\frac{z}{j_{\nu,n}^4}\right).$$
Taking the logarithmic derivative of both sides of the next expression
$$2^{2\nu}z^{-2\nu}\Gamma^2(\nu+1)J_{\nu}({z})I_{\nu}({z})=\prod_{n\geq1}\left(1-\frac{z^4}{j_{\nu,n}^4}\right)$$
and using the recurrence relations $zJ_{\nu}'(z)-\nu J_{\nu}(z)=-zJ_{\nu+1}(z)$ and $zI_{\nu}'(z)-\nu I_{\nu}(z)=zI_{\nu+1}(z)$ it follows that
$$-\frac{1}{4}\left(4\nu+\frac{zI_{\nu+1}(z)}{I_{\nu}(z)}-\frac{zI_{\nu+1}(z)}{I_{\nu}(z)}\right)=\sum_{n\geq1}\frac{z^4}{j_{\nu,n}^4-z^4}.$$
On the other hand, it is well-known that the $\nu\mapsto j_{\nu,n}$ is increasing on $(-1,\infty)$ for each $n\in\mathbb{N}$ fixed, and thus $j_{\nu,n}>{\dots}>j_{\nu,1}\geq j_{-\frac{1}{2},1}=\frac{\pi}{2}>1$ for each $\nu\geq -\frac{1}{2}$ and $n\in\{2,3,\dots\}.$ Moreover, the function
$$\nu\longmapsto \sum_{n\geq1}\frac{1}{j_{\nu,n}^4-1}$$
is decreasing on $(-1,\infty),$ which implies that
$$\sum_{n\geq1}\frac{1}{j_{\nu,n}^4-1}\leq 1$$
if and only if $\nu\geq\nu^{\star},$ where $\nu^{\star}$ is the unique root of the equation
$$\sum_{n\geq1}\frac{1}{j_{\nu,n}^4-1}=1\ \ \Longleftrightarrow\ \ J_{\nu+1}(1)I_{\nu}(1)-J_{\nu}(1)I_{\nu+1}(1)=(\nu+1)J_{\nu}(1)I_{\nu}(1).$$
Thus, applying Lemma \ref{lem1} the proof is complete.
\end{proof}

\end{document}